\documentclass[12pt]{amsart} 

\newtheorem{thm}{Theorem}[section]	
\newtheorem{lem}[thm]{Lemma}		
\newtheorem{cor}[thm]{Corollary}	 

\theoremstyle{definition}       %==== definition type declarations
\newtheorem{dfn}[thm]{Definition}	
	
\newtheorem{prob}[thm]{Problem}		\newtheorem{question}[thm]{Question}

\newtheorem{rem}[thm]{Remark}

\newcommand{\rar}{\rightarrow} 	
	
   	\newcommand{\bfc}{{\mathbb{C}}}
\newcommand{\bfq}{{\mathbb{Q}}}	
\newcommand{\bfz}{{\mathbb{Z}}}	
\newcommand{\bfr}{{\mathbb{R}}}	\newcommand{\bfa}{{\mathbb{A}}}

\newcommand{\ttimes}{\widetilde{\times}}

\newcommand{\setmin}{\,{^{_{\setminus}}}\,}
\newcommand{\sk}{{\mathrm{Sk}}}
\newcommand{\conv}{{\mathrm{Conv}}}
\newcommand{\Q}{{\mathbb{Q}}}
\newcommand{\R}{{\mathbb{R}}}
\newcommand{\Rn}{\R^n}
\newcommand{\C}{{\mathbb{C}}}
\newcommand{\N}{{\mathbb{N}}}
\newcommand{\Z}{{\mathbb{Z}}}
\newcommand{\Zn}{\Z^n}

\newcommand{\Csn}{(\C^*)^n}
\newcommand{\cC}{{\mathcal{C}}}

\newcommand{\bO}{{\mathbf{O}}}
\newcommand{\cone}{{\mathrm{Cone}}}
\newcommand{\cP}{{\mathcal{P}}}
\newcommand{\eps}{{\varepsilon}}

\newcommand{\spec}{{\operatorname{Spec }}}

\newcounter{cumulative}

\begin{document}

\renewcommand{\thefootnote}{\fnsymbol{footnote}}

\noindent%
{\large{\bf
Extending Triangulations and Semistable Reduction\\[2mm]}
\normalsize% 
%\parbox{3.2in}
{D. Abramovich\footnote{Partially supported by NSF grant 
DMS-9503276 and an Alfred P. Sloan research fellowship.}\\ 
\small
Department of Mathematics, Boston University\\
111 Cummington, Boston, MA 02215, USA\\
{\tt abrmovic@math.bu.edu}} \\[2mm]
%\parbox{3.2in}
\normalsize% 
{J. M. Rojas\footnote{Partially supported by an NSF postdoctoral 
fellowship.}\\ 
\small
Department of Mathematics, Massachusetts Institute of Technology\\
 Cambridge, MA 02139, USA \\
{\tt rojas@math.mit.edu}}\\[1mm]}
\noindent %
PRELIMINARY VERSION,  
\today   
%
%\large
\addtocounter{section}{-1}
\addtocounter{footnote}{-2}
\renewcommand{\thefootnote}{\arabic{footnote}}
\section{INTRODUCTION}
In the past three decades, a strong relationship has been established between
convex geometry, represented by convex polyhedra and polyhedral complexes, and
algebraic geometry, represented by toric varieties and 
toroidal embeddings. In this note we exploit this 
relationship in the following manner. We address a basic
problem in algebraic geometry: a certain version of 
semistable reduction. We translate a local  
case of the problem into a basic problem about polyhedral complexes: extending 
triangulations. Once we solve the second problem, the first follows.
We have taken the opportunity with this note to try to extend some bridges
between the terminologies of these two theories.

\subsection{Semistable Reduction} We work over the field of complex numbers
$\bfc$.  Let $f:X\to B$ be a proper morphism of algebraic varieties, whose
generic fiber is reduced and absolutely irreducible. Thus there exists a
Zariski dense open set $U\!\subset\!B$ such that the fiber $f^{-1}(b)$ over any
point in $b\in U$ is a compact complex algebraic variety.

Loosely speaking, semistable reduction for a morphism like $f$ is a
meta-problem of ``desingularization of morphisms,'' where the goal is to
``change $f$ slightly'' so that it becomes  ``as nice as possible''. Of course,
we need to specify more precisely what we mean by the clauses in quotation
marks. 

\subsubsection{What do we mean by a morphism being ``as nice as
possible?''} 

First of all, $X$ and $B$ should be as nice as possible, namely
nonsingular. Moreover, we want $f$ to have a nice, explicit local description,
so that the fibers of $f$ have the simplest possible singularities. 

Such a wonderful morphism will be called {\bf 
semistable}. Here is the definition:

\begin{dfn}
Let $f:X\to B$ be a flat projective morphism, with connected fibers, of 
nonsingular varieties. We say that
$f$ is {\bf semistable} if for each point $x\!\in\!X$ with $f(x)\!=\!b$ there 
is a
choice of formal coordinates $B_b = \spec\ \C[[t_1,\ldots,t_m]]$ and 
$X_x = \spec\ \C[[x_1,\ldots,x_n]]$, such that $f$ is given by: 
$$t_i = \prod_{j=l_{i-1}+1}^{l_i} x_j,$$
where $0 = l_0 < l_1 \cdots < l_m \leq n$, $n=\dim X$, and $m = \dim B$.
\end{dfn}

We must state right up front that in this note we will {\bf not} end up with a
semistable morphism, but we will get very close. In particular, 
our results here form an additional step in recent work on semistable 
reduction \cite{aj,ak}. 

\subsubsection{What do we mean by ``changing $f$ slightly?''} 

First we
define two types of operations necessary for semistable reduction:

\begin{dfn} An {\bf alteration } $B_1\to B$ is a proper, generically finite,
surjective morphism. A {\bf modification} $Y\rar X$ is a birational proper
morphism (equivalently, a birational alteration).
\end{dfn}

Given a morphism $X\to B$ as before, and an alteration $B_1\to B$, we call the
component of 
$X\times_B B_1$ dominating $B_1$ the {\bf main component} and 
denote it by $X\ttimes_B B_1$. 

We are now ready to state the semistable reduction problem in its ultimate
form: 

\begin{prob}\label{semistable-reduction-problem}
Let $X\rar B$ be a flat projective morphism, with 
connected fibers, of nonsingular varieties. Find an
alteration $B_1\rar B$, and a modification 
 $Y\rar X\ttimes_B B_1$, such that $Y\rar B_1$ is semistable.
\end{prob}

\subsubsection{Nearly Semistable Morphisms}
We will need some terminology in order to state the 
weaker version of semistable reduction we actually 
address here. We will follow \cite{te} for the basic 
definitions.\footnote{Also, mimicking standard notation from algebraic
topology,  
$f : (X,A) \longrightarrow (Y,B)$ will be understood to mean that 
$A$ and $B$ are subvarieties of $X$ and $Y$ respectively; and that  
$f$ is a morphism from $X$ to $Y$ satisfying $f(A)\subset B$.} 

\begin{dfn}  \mbox{}\\

\vspace{-.5cm}
\begin{enumerate}
\item A {\bf toric variety} is a normal\footnote{Although 
normality is not assumed in some contexts, all toric 
varieties will be normal in this paper.}  variety 
$X$ with an open embedded copy $T$ of $\Csn$, such that 
the natural $\Csn$-action on $T$ extends to all of $X$. 
We sometimes call the pair $(X,T)$ a {\bf torus 
embedding}. 
\item More generally, suppose $Y$ is a normal 
variety with a smooth open subvariety $U_Y$ 
satisfying the following condition: locally 
analytically at every point, $(Y,U_Y)$ is 
isomorphic to a local analytic neighborhood of 
some torus embedding $(X,T)$. We then call $Y$ 
a {\bf toroidal variety} and $(Y,U_Y)$ a {\bf toroidal 
embedding}.\footnote{We will sometimes 
follow \cite{te} and also refer to the inclusion 
$U_Y\subset Y$ as a toroidal embedding.}
\item A dominant morphism $f: (X,U_X) \to (B,U_B)$ of 
toroidal embeddings is called a {\bf toroidal 
morphism}, if locally analytically near every point on 
$X$ it is isomorphic to a torus equivariant morphism of 
toric varieties.   
\end{enumerate}
\end{dfn}

Roughly speaking, a toric variety is ``monomial:'' an affine toric variety is
always defined by binomial equations, and any toric variety can always be
covered by affine charts in such a way that every overlap isomorphism is  
a monomial map. Similarly, a toroidal variety is ``locally monomial'' and a
toroidal morphism is a ``locally monomial morphism.'' 

If $U_B \subset B$ is a toroidal embedding, then we may write $B\setmin U_B$
as a union of divisors $D_1 \cup\cdots \cup D_k$. More precisely, 
recall that $B\setmin U_B$ can be decomposed into 
strata of varying dimensions (see \cite{te} or \cite{gormac}). In particular,  
let us define $U_B^{(2)}$ to be the union of $U_B$ and the 
codimension 0 strata of $B\setmin U_B$. This notation makes 
sense since we've actually only removed pieces of codimension $\geq\!2$ 
from $B$ to construct $U_B^{(2)}$.

We now detail the type of morphisms we will treat:
\begin{dfn}
A proper toroidal morphism $f: (X,U_X) \to (B,U_B)$ is said to
be {\bf nearly semistable} if the following conditions hold:
\begin{enumerate}
\item There are no horizontal divisors in $X$, namely:  $U_X = f^{-1}(U_B)$.
\item The base $B$ is nonsingular.
\item The morphism $f$ is equidimensional.
\item All the fibers of $f$ are reduced.
\item The restriction of $f$ to $U_B^{(2)}$ is semistable, i.e., ``$f$ is
semistable in codimension $\leq\!1$.'' 
\item The singularities of variety $X$ are at worst finite quotient
singularities.
\end{enumerate}
\end{dfn}

One may ask how far a nearly semistable morphism is 
from a semistable one. The answer is simple: every 
toroidal semistable morphism is nearly semistable; and a nearly semistable
morphism $X 
\to B$ is semistable if and only if $X$ is nonsingular (see \cite{ak}). 
%This is the gist of Theorem \ref{torsemi} in the next section. 
% Or is this ``simple answer'' quoted from some 
% reference?  If so, what's the reference?

\subsubsection{The Result}

The problem addressed in this paper is a special (local) case of nearly
semistable reduction: 

\begin{thm}\label{torsemi}
Set $B\!=\!\bfa^n_\C$ and let $U_B$ be the natural open 
subscheme of $B$ whose underlying complex variety is 
$\Csn$. Note that the 
inclusion $U_B \subset B$ is a toroidal embedding, and 
let $f:X\to B$ be a proper morphism satisfying:
\begin{enumerate} 
\item $U_X:=f^{-1}(U_B) \subset X$ is a toroidal 
embedding,
and $f:(X,U_X) \to (B,U_B)$ is a toroidal morphism;
\item $f$ is equidimensional, with smooth and absolutely irreducible generic
      fiber; 
\item every fiber of $f$ is reduced.
\end{enumerate} 
Then there exists a {\bf finite} toric morphism $(B_1,U_{B_1})\to
(B,U_B)$ and a toroidal modification $Y\to X\times_B B_1$, such that 
$Y\to B_1$ is nearly semistable.
\end{thm}

One may ask what right we have to make all these assumptions on the
morphism $f$ we start with. In \cite{ak} it is shown that given any morphism
$f$, as in Problem \ref{semistable-reduction-problem}, we can reduce it to a 
toroidal morphism $f$ as in Theorem \ref{torsemi}. Such morphisms are 
called {\bf weakly semistable} in \cite{ak}.

The methods of \cite{ak} are quite different from what we do here. In short,
they involve: 
\begin{enumerate}
\item Making $X\to B$ toroidal. This follows easily from the methods of
\cite{aj}. 
\item Making a toroidal $X\to B$ satisfy the conditions 
in the theorem. Locally this can be done easily using 
toroidal modifications and finite base changes. To do 
it globally one uses a covering trick of Kawamata (see
\cite[Theorem 17]{kawamata}). 
\end{enumerate}
Moreover, once the local results here are established, we can go back to 
\cite{ak} and, using Kawamata's covering trick, extend it to prove nearly
semistable reduction in general. 

\subsection{Extending Triangulations} 

We now wear our polyhedral glasses. 

For the concepts of a {\bf compact polyhedral complex} 
$\Delta$ and a {\bf conical polyhedral complex} 
$\Sigma$ see \cite[pg.\ 69, Definition 5]{te}.  
An {\bf integral structure} on a compact or conical 
polyhedral complex is defined in \cite[pg.\ 70, 
Definition 6]{te}. We will always assume that our 
complexes come equipped with an integral structure. 
{}From here on, we will simply say {\bf polyhedral 
complex}, when we mean a compact polyhedral complex 
with integral structure. 

\begin{rem} 
A useful example of a polyhedral complex to consider is a finite collection  
$\cP$ of integral polyhedra in $\Rn$. (Recall 
that a polyhedron in $\Rn$ is {\bf integral} iff 
all its vertices lie in $\Zn$.) If  
$\cP$ is closed under intersection and taking 
faces, then $\cP$ is a polyhedral complex. 
Note, however, that {\bf not} all polyhedral 
complexes arise this way. This accounts 
for some of the geometric richness of 
toroidal varieties.
\end{rem}

Again, in \cite[pg.\ 70]{te}, it is shown that for any 
compact polyhedral complex $\Delta$, one can construct 
a conical polyhedral complex, which we denote 
$\Sigma(\Delta)$ --- namely the cone over $\Delta$. To 
reverse the process, define a {\bf slicing function} 
$h: \Sigma \to \bfr$ to be a nonnegative continuous 
function, whose restriction to every cone $\sigma\in 
\Sigma$ is linear, which vanishes only at
the origin $\bO\!\in\!\Sigma$. Then the {\bf slice} 
$h^{-1}(1)$ of $\Sigma$ defines a compact polyhedral 
complex $\Delta(\Sigma,h)$.

We denote by $\sk^k(\Delta)$ the $k$-skeleton of 
$\Delta$. We will also use $\#S$ for the cardinality of 
a set $S$, and $\cone(V)$ for the set of all 
nonnegative linear combinations of a set of 
vectors $V\subset\Rn$. 

By a  {\bf subdivision} $\Delta'$ of $\Delta$ (resp.\ 
$\Sigma'$ of $\Sigma$) we will
mean a finite partial polyhedral decomposition 
of $\Delta$ (resp.\ $\Sigma$), as in \cite[pg.\ 86, 
Definition 2]{te}, with the {\bf 
completeness} property: $|\Delta'| = |\Delta|$. 
(Recall that the notation $|\Delta|$ simply 
means the topological space consisting of the 
union of all the cells of $\Delta$.)
A subdivision $\Delta'$ is called a {\bf triangulation} 
or a {\bf simplicial subdivision} if every cell of 
$\Delta'$ is a simplex.

A {\bf lifting function} (or {\bf order function}) 
$f: \Delta \to \bfr$ on a polyhedral complex 
is a continuous function, convex and piecewise linear 
on each cell of $\Delta$, respecting the integral 
structure. In the conical case ($f: \Sigma \to \bfr$) 
we add the requirement that $f$ be homogeneous: 
$f(\lambda x)=\lambda f(x)$, for all $\lambda\!\geq0$ 
and all $x\!\in\!|\Delta|$ \cite[pg.\ 91, condition 
($*$)]{te}. 

\begin{rem}\label{convexity} We follow the convention in \cite{te}, where  one
requires a lifting  
function to be ``convex down'' on each cell, namely $f(\lambda x + \mu y) \geq
\lambda f(x) + \mu f(y)$.  Also, all our lifting functions take rational
values on the lattices in the cells. 
This is in contrast with the polyhedral convention, as in \cite{z}, 
where lifting functions are ``convex up'' and real values are allowed. 
\end{rem}

Given a lifting function $f: \Delta \to \bfr$, (resp.\ 
$f:  \Sigma \to \bfr$) we define the subdivision 
$\Delta_f$ (resp.\ $\Sigma_f$) {\bf induced by $f$}, to 
be the coarsest subdivision such that $f$ is linear on 
each cell.\footnote{Here we are poised on the verge of 
a notational quagmire. In \cite{te}, a subdivision 
induced by a lifting function is called {\bf 
projective}, which makes sense from the 
algebro-geometric point of view, as the corresponding  
modification of toric varieties $X_{\Sigma'} \to X_\Sigma$ is a projective 
morphism. However, this kind of projectivity is foreign to the polyhedral 
world. In \cite{z} such subdivisions are called {\bf 
regular}, which is fine by us, except that in \cite{te} 
the term ``regular'' is sometimes used for an explicit 
type of subdivision which is very regular indeed... We 
will thus simply refer to our subdivisions as ``subdivisions induced 
by lifting functions.''}

\begin{rem}\label{induced}
The subdivision induced by $f$ is clearly determined by 
the values of $f$ on its vertices $\sk^0(\Delta_f)$ 
(resp.\ its edges $\sk^1(\Sigma_f)$). In fact one can construct $f$ from its
values on $\sk^0(\Delta_f)$ (resp.  $\sk^1(\Sigma_f)$) as the minimal function
which is convex on each cell, having the given values on  $\sk^0(\Delta_f)$
(resp.  $\sk^1(\Sigma_f)$).  
 In particular, a subdivision induced by a lifting function 
{\bf can} sometimes add new vertices to $\Delta$ 
(resp.\ new edges to $\Sigma$). However, with 
some care, we can control this behavior.
\end{rem}

We will prove the following result:
\begin{thm}\label{regular-extension}
Let $\Delta$ be a polyhedral complex and $\Delta_0\subset \Delta$
a subcomplex. Let $\Delta_0'$ be a triangulation of $\Delta_0$ induced by a
lifting function. Then there exists a triangulation 
$\Delta'$ of $\Delta$, also induced by a lifting 
function, which extends $\Delta_0'$ and introduces 
no new vertices. That is, $\sk^0(\Delta')=
\sk^0(\Delta)\cup\sk^0(\Delta'_0)$.  
\end{thm}

Applying this to a slice of a conical polyhedral complex we obtain:
\begin{cor}\label{conical-regular-extension}
Let $\Sigma$ be a conical polyhedral complex admitting 
a slicing function $h:\Sigma\to \bfr$, and let 
$\Sigma_0\!\subset\!\Sigma$ be a subcomplex. Let 
$\Sigma_0'$ be a triangulation of $\Sigma_0$ induced by 
a lifting function. Then there exists a triangulation 
$\Sigma'$ of $\Sigma$, also induced by a lifting
function, which extends $\Sigma_0'$ and introduces 
no new edges.  That is, $\sk^1(\Sigma')=\sk^1(\Sigma)
\cup \sk^1(\Sigma'_0)$. 
\end{cor}

One may ask, ``Do we really need to assume that $\Delta'_0$ is induced by a 
lifting function?'' The simplest example showing that 
this is indeed the case was communicated to us 
independently by R. Adin and B. Sturmfels:

Let $\Delta\subset \bfr^3$ be the triangular prism $\delta =
\conv\{f_{0,0},\ldots,f_{1,2}\}$, where:
$$\begin{array}{rclrclrcl}
   f_{0,0} & = & (0,0,0);\quad &
   f_{0,1} & = & (1,0,0);\quad &
   f_{0,2} & = & (0,1,0) \\
   f_{1,0} & = & (0,0,1);\quad &
   f_{1,1} & = & (1,0,1);\quad &
   f_{1,2} & = & (0,1,1) 
\end{array} $$

Let $\Delta_0 = \partial \Delta$ be the boundary of our prism.

Let  $\Delta_0'$ be the subdivision of $\Delta_0$ obtained by inserting 
the following new edges: $$\overline{f_{0,0}f_{1,1}},
\overline{f_{0,1} f_{1,2}},\overline{f_{0,2}f_{1,0}}$$ (So we've 
``cut'' a new edge into each $2$-face of $\Delta_0$.) It is an easy exercise 
to see that there is no extension of $\Delta_0'$ (to a triangulation of 
$\Delta$) without new vertices: in particular, any $3$-cell of such an 
extension must have an edge intersecting the midpoint of some edge of 
$\Delta'_0$ --- a contradiction. It is also not hard to 
see that $\Delta_0'$ can not be induced by any lifting 
function \cite[Chapter 3]{fulton}.

\section{Reduction of Theorem \ref{torsemi} to \ref{regular-extension} 
%or \ref{fiber-extension} 
}

Let $f:X\to B$ be as in  Theorem \ref{torsemi} and $f_{\Sigma}:\Sigma_X
\to \Sigma_B$ the associated morphism of rational 
conical polyhedral 
complexes. Note that $\Sigma_B$ is a nonsingular cone (a simplicial cone of 
index 1): it is simply  
the nonnegative orthant in $\bfr^n$, generated by the 
standard basis vectors $\{\hat{e}_i\}$.  Let $\tau_i$ be
the edges of $\Sigma_B$, namely $\tau_i = 
\cone(\hat{e}_i)$.  We identify the lattice of $\tau_i$ 
with $\bfz \hat{e}_i$. 

Let $\Sigma_B^1 = \bigcup \tau_i$ be the 1-skeleton of $\Sigma_B$ and 
$\Sigma_X^1 = f_\Sigma^{-1}(\Sigma_B^1)$.  Also let $\Sigma_{X,i} =
f_\Sigma^{-1}(\tau_i)$.  For an integer $k_i$ let $N_i(k_i)$ be the integral 
structure on $\Sigma_{X,i}$ obtained by intersecting the lattices in 
$\Sigma_{X,i}$ with $f_\Sigma^{-1}(\bfz k_i\cdot \hat{e}_i)$.

By \cite[Chapter III, Theorem 4.1 pg.\ 161]{te}, as 
interpreted in \cite[Chapter II, \S 3]{te}, there exists 
an integer $k_i$ and a simplicial subdivision 
$\Sigma'_{X,i}$ of $\Sigma_{X,i}$, {\bf which is induced by a 
lifting function}, having index 1 with respect to the 
integral structure $N_i(k_i)$.

Let $B_1\simeq \bfa^n_\C$ be complex affine space with 
coordinates $s_1,\ldots, s_n$. The substitution 
$s_i^{k_i} = t_i$ gives a homomorphism 
$\bfc[t_1,\ldots,t_n] \to \bfc[s_1,\ldots,s_n]$, giving 
rise to a finite morphism  $B_1\to B$. Then 
$\Sigma_{B_1}$ is the same as $\Sigma_B$ 
but taken instead with the lattice
$N_{B_1}=\prod \bfz k_i\hat{e}_i$. Let $X_1 = 
X\times_B B_1$. Since the fibers of $X$ are reduced, it follows that
$X_1$ is normal and $X_1\to B_1 $ is again toroidal.
Likewise, $\Sigma_{X_1}$ is just $\Sigma_X$ with 
integral structure given by intersecting the lattices 
in $\Sigma_X$ with $f_\Sigma^{-1}(N_{B_1})$. 

Putting the triangulations $\Sigma'_{X,i}$ of 
$\Sigma_{X,i}$ together, there exists a triangulation 
$\Sigma_X^{1'}$ of $\Sigma_X^1$ (induced by a lifting 
function) of index 1 with respect to the integral 
structure on $\Sigma_{X_1}$!

Let us verify that $\Sigma_X$ admits a slicing function:
let $h_B:\Sigma_B \to \bfr$ be the function defined by 
$h_b(\sum a_i \hat{e}_i)=  \sum a_i$. Then the pullback 
$h_b\circ f_\Sigma$ is a slicing function on
$\Sigma_X$.

By Corollary \ref{conical-regular-extension} of Theorem
\ref{regular-extension}, there is an extension of 
$\Sigma_X^{1'}$ to a triangulation $\Sigma_X'$ of 
$\Sigma$ (induced by a lifting function) without added 
edges. 

Let $Y\to X_1$ be the corresponding toroidal 
modification and let $f_1 : Y\to B_1$ the resulting 
morphism.

Note that since all the edges in the triangulation $\Sigma_X'$ map to the edges
$\tau_i$ of $\Sigma_{B_1}$, we have that $f_1$ is equidimensional
\cite{ak}. Since
the integral generator of every edge in $\Sigma_X'$ maps to the generator of
the image edge in $\Sigma_{B_1}$, and since $B$ is nonsingular, all the fibers
of $f_1$ are reduced \cite{ak}. By the construction of \cite{te}, $f_1$ is 
semistable in codimension 1. Since $\Delta_X'$ is simplicial, $Y$ has at worst 
quotient singularities. Thus $f_1$ is nearly semistable. \qed

\begin{rem}
The variety $Y$ may be singular, as the following 
example shows: let $\Sigma_Y\subset \bfr^4$ be the 
nonnegative orthant, generated by the standard basis 
vectors $\hat{e}_1,\ldots \hat{e}_4$. Let 
$w\!=\!(1/2,1/2,1/2,1/2) \in  \bfr^4$ and $N_Y$ the 
lattice generated by $w,\hat{e}_1,\ldots \hat{e}_4$. 
Also let $Y$ be the corresponding toric variety --- 
the quotient of $\bfa^4_\bfc$ by the diagonal 
$\bfz/2$ action given by $p \mapsto -p$ --- which 
happens to be singular. Finally, let 
$\Sigma_B\subset \bfr^2$ be the first quadrant, 
generated by the standard basis vectors $\hat{e}_1, 
\hat{e}_2$, with the standard lattice 
$N_B=(\{0\}\cup\N)^2$. We have a canonical morphism 
$\Sigma_Y \to \Sigma_B$ via $$(a,b,c,d) \mapsto 
(a+b,c+d)$$ which maps $N_Y$ into $N_B$. The resulting 
morphism $Y \to \bfa^2_\bfc$ is nearly semistable, but 
not semistable.
\end{rem}

\section{Proof of Theorem \ref{regular-extension}}

%Goal:
%Given a complex $\Delta$, a subcomplex $\Delta'$ and a regular triangulation
%$\Delta'_0$ of 
%$\Delta_0$ (possibly with added vertices), then there is a regular
%triangulation 
%$\Delta'$ without additional vertices extending $\Delta'_0$. 

% We can of course assume $\Delta_0\!\neq\!\Delta$, for 
% otherwise we would already be done.
% Also, we may assume 
% that $\Delta$ is connected, by the definition of an 
% induced subdivision. 

It is a simple fact, made precise in Lemma \ref{lem:gen} below, 
that any {\bf generic} lifting function on a polyhedral 
complex induces a simplicial subdivision. This fact 
is used frequently in applications of subdivisions to 
the computation of mixed volumes, polyhedral homotopies, 
and {\bf toric} (or {\bf sparse}) resultants 
\cite{sparseelim,polyhomo,ce,affres}. The last 
two constructions give effective recent techniques, 
sometimes  more efficient than Gr\"obner bases, 
for solving systems of polynomial equations. 

However, it should be emphasized that the lifting 
functions considered here and in \cite{te} are more 
general than those in 
\cite{sparseelim,polyhomo,ce,affres}: 
via the use of convex hulls, the lifting 
functions there are completely determined by the values 
assigned to the vertices of $\Delta$. We will call 
these more restricted lifting functions {\bf 
verticial}. The verticial lifting functions are a bit 
more ``economical'' in the sense that their 
corresponding subdivisions never introduce any new 
vertices.  

There is an easy way to resolve this difference by passing to the 
verticial case from the start. In fact, we will reduce the 
proof of Theorem \ref{regular-extension} to finding {\bf any} triangulation 
(given by a verticial lifting function) 
in a new, specially constructed, polyhedral complex. The 
latter problem is then almost trivial to solve. 

First recall (see \cite{te}, Corollary 1.12) that 
induced subdivisions are transitive:  if $\Delta'$ is a subdivision of 
$\Delta$ induced by a lifting
function $f$ on $\Delta$, and $\Delta''$ is a subdivision of $\Delta'$ induced
by a lifting function $f'$ on $\Delta'$, then $\Delta''$ is a subdivision of
$\Delta$ as well. In fact, $\Delta''$ is induced by $f + \epsilon f'$ for 
sufficiently small $\epsilon>0$. 

Thus let $f_0:\Delta_0\to \bfr$ be a lifting function which induces
the given subdivision $\Delta_0'$ in our theorem. By adding a
constant if necessary, we may assume $f_0$ is positive. Following Remark 
\ref{induced}, we can take the values of $f_0$ on $\sk^0(\Delta_0')$, extend
them by zero to the other vertices $\sk^0(\Delta)\setmin \sk^0(\Delta_0')$, and
take the minimal lifting function $f:\Delta\to \bfr$ which has these values on
the vertices $\sk^0(\Delta)\cup \sk^0(\Delta_0')$. Clearly $f|_{\Delta_0} =
f_0$. Let $\Delta_1$ be the induced subdivision. Then clearly the restriction
of $\Delta_1$ to $\Delta_0$ coincides with $\Delta'_0$. If $\Delta'$ is any
subdivision of $\Delta_1$ without new vertices, then its restriction 
to $\Delta_0$ must be $\Delta'_0$, since $\Delta_0'$ is already simplicial: any
subdivision of a simplicial complex without new vertices is trivial. Thus all
we need to do to prove Theorem \ref{regular-extension} is find a {\bf 
verticial} lifting function 
on $\Delta_1$ giving a triangulation. In summary, by replacing $\Delta$ with 
$\Delta_1$, we can assume that $\Delta_0=\Delta_0'$ and then  
conclude by finding {\bf any} triangulation of $\Delta_1$ 
(given by a verticial lifting function) --- 
a simpler problem than finding a triangulation of one complex 
extending some other triangulation. 

To complete the proof of Theorem \ref{regular-extension}, recall the following 
lemma:

\begin{lem}
\label{lem:gen}
Supppose $\Delta$ is a polyhedral complex. Then 
\begin{enumerate}
\item{The set $L_\Delta$ of all verticial lifting  functions on 
$\Delta$ is a finite-dimensional rational vector space.} 
\item{The set of all lifting functions which do 
{\bf not} induce simplicial subdivisions is a finite 
union of proper subspaces of $L_\Delta$. }
\end{enumerate}
\end{lem}

\noindent
{\bf Proof:} Note that any verticial lifting function on 
$\Delta$ is uniquely determined by its values on 
$\sk^0(\Delta)$, which are assumed to be rational, so part (1) 
follows immediately.

To prove (2), let $\cC:=(c_v \; | \; v \in \sk^0(\Delta) )$ be a vector of
rational constants.  Let $\Delta_\cC$ denote the subdivision of $\Delta$
induced by the verticial lifting function sending $v \mapsto c_v$ for  
all $v\!\in\!\sk^0(\Delta)$. 

Now suppose that there is a nonsimplicial cell $C$, 
with vertex set $V(C)$, in $\Delta_\cC$. Recall that 
the coordinates of $d+2$ points lying on a $d$-flat 
in $\Rn$ must satisfy a determinant depending only on 
$(d,n)$.\footnote{For example, $(x_1,y_1)$, 
$(x_2,y_2)$, and $(x_3,y_3)$ lie on a line 
iff $\begin{vmatrix} x_2-x_1 & y_2-y_1 \\
                     x_3-x_1 & y_3-y_1 
\end{vmatrix}=0$. Note also that this determinant 
is linear in the ``last'' coordinates $\{y_i\}$.} 
(In particular, this determinant is a nonconstant 
multilinear function in the coordinates of the 
points.) Then, by the definition of a cell in a 
subdivision induced by lifting, there must be a 
(nontrivial) linear relation satsified by $(c_v \; | 
\; v\!\in\!V(C))$. 
Furthermore, this linear relation depends only on 
$\Delta$ and the set of vertices $V(C)$. Since 
there are only finitely many possible nonsimplicial 
cells (since, by definition, our polyhedral complexes have only finitely 
many vertices), (2) follows immediately. \qed 

The following is an immediate corollary of our lemma. 
\begin{cor}
\label{cor:good}
Recall the notation of the proof of Lemma \ref{lem:gen}, 
and endow $\bfq^{\#\sk^0(\Delta)}$ with the standard 
Euclidean metric $\|\cdot\|$. Let $\cC\in\Q^{\#\sk^0(\Delta)}$. 
Then for sufficiently small $\eps>0$, 
\begin{enumerate}
\item{$\Delta_{\cC'}$ is a simplicial subdivision 
for {\bf some} $\cC'\!\in\!\bfq^{\#\sk^0(\Delta)}$ satisfying 
$\|\cC'-\cC\|<\eps$. }
\item{If $\Delta_{\cC}$ is already a simplicial 
subdivision, then so is $\Delta_{\cC'}$, {\bf for all} 
$\cC'\!\in\!\bfq^{\#\sk^0(\Delta)}$ satisfying 
$\|\cC'-\cC\|\!<\!\eps$.}
\end{enumerate}
\qed 
\end{cor}
\begin{rem}
Put another way, simplicial subdivisions are a  
dense (via (1)) and open (via (2)) subset of  
the space of all subdivisions arising from verticial
lifting functions.  In fact, we really have 
the stronger statement that the set of all 
lifting values giving a {\bf particular} simplicial 
subdivision forms an open cell within the 
space of all subdivisions.  

Note also two ``nearby'' subdivisions $S_1$ and 
$S_2$ need not have the same extensions, even 
if $S_1=S_2$: for example, consider the unit 
square $S$ with vector of vertices (ordered 
clockwise) $(a,b,c,d)$, and the subcomplex $E$ 
consisting of the edges $\{a,b\}$ and $\{c,d\}$. 
Then $\cC=(0,0,0,0)$ and $\cC'=(-1,1,-1,1)$ 
both generate the same (trivial) subdivision
of $E$. However, these two liftings generate 
different subdivisions of $S$, the first 
being trivial.  
\end{rem}

Returning to the proof of Theorem 
\ref{regular-extension}, it follows by Corollary \ref{cor:good}
that there exists a simplicial subdivision of $\Delta_1$ without new vertices,
which is what we needed to prove. \qed 
% let $\{\Delta_0(i)\}$ be the (finitely many) connected 
% components of $\Delta_0$. Also let 
% $\cC'_0$ be the 
% vector of lifting values inducing $\Delta'_0$. 
% Note that by Corollary \ref{cor:good}, any small
% perturbation of $\cC'_0$ will still give the 
% same subdivision $\Delta'_0$. 

% Let us now alter our complexes slightly: 
% Redefine $\Delta$ (using the coarsest possible 
% subdivision) so that every vertex of $\Delta'_0$ 
% is a vertex of $\Delta$. (In particular, 
% every vertex of $\Delta'_0$ is now also a vertex of 
% $\Delta_0$.) We may now work solely with verticial 
% lifting functions.  

%Define $\cC'$ to be any point in 
%$\R^{\#\sk^0(\Delta)}$ such 
%that the subvector corresponding to $\sk^0(\Delta_0)$ is 
%identical to $\cC'_0$. Note that by the definition of an 
%induced subdivision, $\Delta'_{\cC'}|_{\Delta_0} 
%=\Delta'_0$. So $\Delta'_{\cC'}$ is already a bona fide 
%extension of $\Delta'_0$. The only trouble 
%is making it satisfy the conclusion of our theorem.   

%Now if $\Delta'_{\cC'}$ is a simplicial subdivision, 
%then we are done. Otherwise, by Corollary 
%\ref{cor:good}, there is some perturbation $\cC''$ 
%of $\cC'$ for which $\Delta_{\cC''}$ is a 
%simplicial subdivision. Similarly, {\bf any} 
%sufficiently small perturbation $\cC''$ of $\cC'$ must 
%satisfy $\Delta_{\cC''}|_{\Delta_0(i)}=\Delta'_0(i)$ 
%for all $i$. Therefore there must be a perturbation 
%$\cC''$ with $\Delta_{\cC''}$ satisfying 
%the conclusion of Theorem \ref{regular-extension}. So 
%we are done. \qed 

\subsection{Acknowledgements} 

We would like to thank B. Sturmfels for inciting this 
collaboration, and R. Adin for a discussion of ideas 
relevant to this note.

%\section{Proof of Theorem \ref{fiber-extension}} (IF IT IS TRUE) 

\end{document}